\documentclass{fundam}

\usepackage{graphicx}
\usepackage{tikz}
\usepackage{amsmath,amssymb}
\usepackage{url}
\usepackage[T2A,T1]{fontenc}
\usepackage[utf8]{inputenc}
\usepackage[russian,english]{babel}
\usepackage{color}

\begin{document}

\setcounter{page}{219}
\publyear{22}
\papernumber{2127}
\volume{186}
\issue{1-4}

   \finalVersionForARXIV

\title{Boris (Boaz) Trakhtenbrot --- The Beginning}

\author{Mark Trakhtenbrot\thanks{Address for correspondence:  Department of Computer Science,  Holon Institute of Technology,
                                   52 Golomb st., Holon 5810201,  Israel}
\\
 Department of Computer Science\\
 Holon Institute of Technology \\
 52 Golomb st., Holon 5810201,  Israel \\
  markt@hit.ac.il
  }

\maketitle

\runninghead{M. Trakhtenbrot}{Boris (Boaz) Trakhtenbrot --- The Beginning}

\vspace*{-12mm}
\section{Introduction}

20 February, 2021 was the centenary date of Boris Abramovich (Boaz) Trakhtenbrot. As it was with most people of his generation, he was living his life on the background of many events of which the 20th century was rich with, and that directly affected him and seriously influenced his life. Among them were global historical events, as well as those that specifically related to the fate of his family and directly to his fate, both personal and scientific. The life path of B.T. to his universally recognized scientific achievements passed through different periods - joyful, dramatic, and even tragic. The goal of this article is to tell about this, with emphasis on the little-known aspects of B.T.'s life related to his roots, background and the beginning of his rise as a scientist. This article covers the first third of his life - from the birth in Bessarabia in 1921 until 1951, when after obtaining his PhD degree he and his young family arrived to town Penza in Central Russia.

To begin with, in the long life of B.T. (he lived for 95 years) there was a variety of changes, sometimes unusual.

\medskip
\textbf{The name changed} --- he was born as Boaz Trachtenbroit, and in 1940 he became Boris Trakhtenbrot. Since then, this name was kept until the end of his life; but after moving to Israel, he gladly added his original name Boaz - and so he became Boris (Boaz) Trakhtenbrot. To the rest of this article, I will use the name Boaz.

\medskip
\textbf{Geography and citizenship changed} --- Romanian, Soviet, Israeli. He was born in Bessarabia, which was part of Romania, and in June 1940 (following the secret protocol of the Molotov-Ribbentrop Pact from August 23, 1939) was annexed to the USSR. And at the end of December 1980,
Boaz made his old dream come true and repatriated to Israel.

\medskip
\textbf{The main languages also changed}. Until 1940, these were Yiddish (native language), Hebrew and Romanian (state language). Since 1940, Russian became the main, and in fact the new native language. After moving to Israel, Hebrew came to the fore. Boaz knew it very well from school age, and he was always proud of this; so he plunged into it with pleasure.

\medskip
\textbf{The areas of his research were constantly changing}, which reflected the breadth of Boaz's scientific interests - descriptive set theory, mathematical logic, automata theory, complexity of computations, frequency computations, semantics of programming languages, concurrency, timed automata, hybrid systems, and more

\smallskip
And now -- to details.\vspace*{-2mm}

\section{Brichevo}\vspace*{-1mm}

Boaz was born in the Jewish settlement (\textit{shtetl}) Brichevo. For him, it was always something much more than just a birthplace. All his life he remembered Brichevo and the Brichevans, various events and people connected with it, the first years of his studies and teachers, and even a popular local witty man, his jokes and expressions. In our family, Boaz used to tell various stories related to Brichevo; we knew them all, as well as the names of people associated with them. It did not stop with Boaz's sons; he also told to his first granddaughter about the life of the Brichevo boy Boaz - and she listened to and absorbed all this with interest. For Boaz, Brichevo was a whole world.

\begin{figure}[!h]
\vspace*{-2mm}
\centering
\includegraphics[scale=0.9]{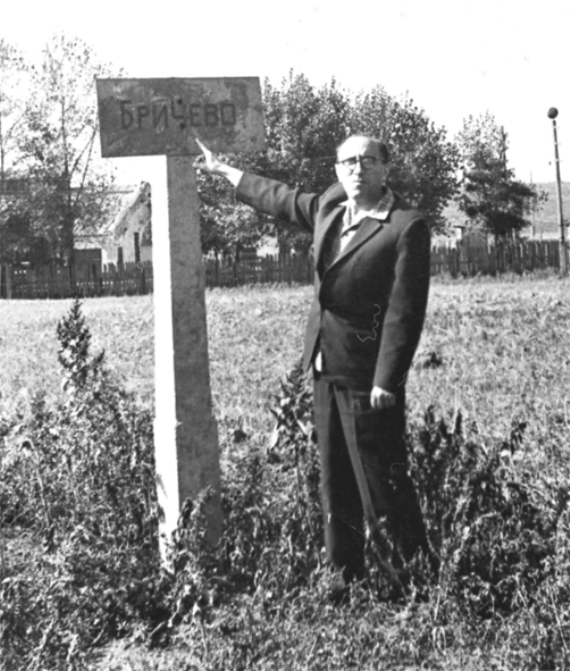}\vspace*{1mm}\\
Boaz visiting Brichevo in 1965
\label{fig:discussion}
\end{figure}

Brichevo was founded in 1836 by a small number of families as a purely Jewish settlement (``colony''). To some extent, it was isolated - villages with a Christian population were at some distance from it. This, according to Boaz's recollections, distinguished Brichevo from other Jewish settlements. Bessarabia was then part of the Russian Empire, to which it was included in 1812 (until 1918, when it became part of Romania). This territory was relatively new for Russia, and the initiative for its development came from the authorities. In particular, there was an idea to give some land to Jewish population. Moreover, the authorities even offered all sorts of benefits to new settlers: tax rebates, reductions in military service (or even exemption from it).

Initially, Brichevo was an agricultural settlement: they grew tobacco, bred small livestock (sheep, goats), made cheese, etc. With time, agriculture-related activities have shrank. Instead, people became involved in various crafts (tailors, shoemakers and blacksmiths), merchant and small agriculture-based industries (oil churns, mills, soap production). They served all the surrounding villages.

By 1897 the population of Brichevo was 1667 inhabitants, and reached a peak in 1930, with 3500 inhabitants; some Ukrainian and Moldovan people joined the population of Brichevo. About that time, economic degradation began, and young people started to look for a better life elsewhere, mainly in South America (in particular, several of Boaz's cousins emigrated to Venezuela). Thus, many families became dependent on the support of their emigrated children. In addition, more ideologically oriented people emigrated to Palestine (the Land of Israel). Following the Holocaust, and later the USSR collapse, the population of Brichevo shrank drastically, with only 305 inhabitants in 2004. 	

\section{Family of Boaz}

Among the founders of Brichevo was Boaz's great-great-grandfather from the side of his mother, Pearl Gelman.
And from the Trachtenbroit side, the first one to come to Brichevo was Boaz's grandfather. They all came from Ukraine. Boaz's grandmother from the Trachtenbroit side was the first one in the family to become involved in merchant activity. Boaz's father, Abraham Trachtenbroit, was the owner of a textile shop; this was when Bessarabia was already a part of Romania. He used to go to a big city (Ia\k{s}i), to buy fabrics there, and then sold it in his shop to the local people of Brichevo. More distant ancestors of Boaz were religious people; but his parents, although they observed all the traditions, were no longer essentially religious.

\begin{figure}[!h]
\centering
\includegraphics[scale=0.9]{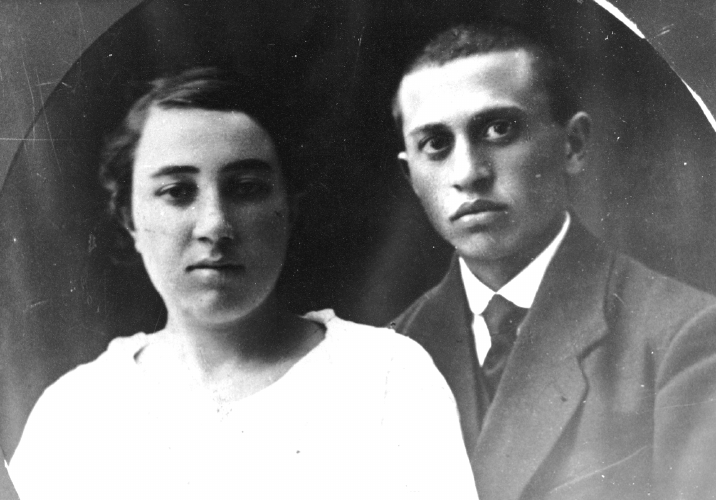}\vspace*{2mm}\\
{Pearl and Abraham Trachtenbroit -- parents of Boaz}
\end{figure}

The family was wealthy; they owned a good house with a garden. But, as in all other homes of Brichevo, there was no electricity and no running water. Abraham permanently donated money for community needs and for help to poor people -- especially because services provided by the state were very restrictive. Therefore, communities had to take care of their needs by themselves.

Boaz's mother Pearl died shortly after she gave birth to Boaz's sister Ratza, when Boaz was three years old. His father then got married for the second time. In this marriage, Boaz's younger brother Mendel was born (in the USSR he used the name Misha). Boaz's stepmother Brucha was a wonderful woman, and she always treated Boaz and his sister as if she was their biological mother. Throughout his life, Boaz used to say that he had two mothers: mother Pearl and mother Brucha.

\section{Cultural aspects of life in Brichevo}

In Brichevo there was a vibrant cultural life, thanks to a group of highly educated people, including autodidacts. Some of these people arrived when they escaped from pogroms in other areas of Russian Empire; for example, pogrom victims fled from Ukraine to Bessarabia (and in particular, to Brichevo) -- there were no pogroms in Bessarabian Jewish settlements.

\medskip
In Brichevo there was a rich library that contained 1800 volumes -- fiction, history books, etc. in Romanian, Yiddish and Hebrew; books in Russian were removed after transition of Bessarabia from Russia to Romania in 1918. There were even meetings with literary and philosophical discussions in the library.

Boaz's father Abraham was a self-educated person; he regularly subscribed to books and newspapers in Yiddish and Hebrew. While Abraham used to read in Yiddish, Boaz's reading was in Hebrew only -- from short stories to books by classic authors (Jules Verne, Arthur Conan Doyle, Leo Tolstoy, etc.), according to his age. But in Brichevo there were also people who could not afford to buy books (they were quite expensive). And Boaz's family had its ``clients''; for example, a neighbor carpenter used to come for taking books and newspapers which the family willingly shared with him.

Theaters were extremely popular in Jewish settlements like Brichevo, and even theaters from Vilna, Warsaw and Chernowitz (including well-known actors such as Sidi Tal) came to perform there. Cantors also were invited to appear. And there were amateur theater groups and musicians that performed for the local public and organized mask-carnivals. Boaz especially remembered a play based on the biblical story of Josef and his brothers.

Parents of Boaz wanted him to play the violin (not very surprising for a boy in a Jewish family), but nothing came out of it. Another attempt was the mandolin, with the same result. He then started to play the flute and really liked it. When he studied in gymnasia, he also played the trumpet in the school orchestra. With this orchestra, he even participated in the 1938 celebrations in Bucharest, on the coronation day of the king of Romania. And many years later, on his 50th birthday in 1971, Boaz got a trumpet as a humorous gift from his former PhD students and colleagues in his department.

\section{School years}

Boaz started his education at the age of five in a \textit{cheder} (``room'' in Hebrew) -- a traditional~elementary school~where boys learned the Hebrew alphabet, reading in the mother tongue (Yiddish), basics of Judaism and prayers in Hebrew. As for prayers -- boys learned them by heart, without understanding the meaning; but this still allowed them to know some Hebrew words. For the lessons, a room was rented in the house of a man who worked in public bathes. As a reward for good achievements, the teacher could pinch a boy's cheek; and a punishment was a hit with a ruler on a boy's hand.

\medskip
After \textit{cheder}, Boaz started his studies in a state (Romanian) elementary school; they continued for four years. In the afternoon, he had classes with a private teacher -- the family could afford to give its children a broader education. These studies were conducted in Hebrew and they included the language and the grammar, literature, history of the Jews, and geography of the Land of Israel. He also studied \textit{Torah}~-- the first five books of the Bible, with the emphasis mainly on biblical stories and not on religious aspects of the Holy Book. Many years later, Boaz told his sons these stories: about Cain and Abel, Jacob and Esau, and his favorite one -- the story of Josef and his brothers. Boaz very much liked all these studies, and throughout his life he remembered with gratitude the teacher, Pinchas Golergant, that taught all these disciplines. Boaz always called him ``my first teacher'', clearly distinguishing him from other teachers -- those that taught him in the state elementary school.

\medskip
When Boaz was in his third year in the elementary school, he got an award. Usually awards were granted yearly for accurate writing and these kinds of skills. Boaz's handwriting was never a candidate for any prize; but his teachers realized that he was a gifted boy. Moreover, the school's director not only praised him very much (according to Boaz -- too much), but also expressed a fear that despite Boaz's excellence in his studies his parents would just send him to work in a shop. Maybe a source of this fear was that Boaz's father wanted him to get an education related to textile technology and to open his own business; in the \textit{shtetl} these plans were definitely known.

Then the time came to continue studies in secondary school --- but there was no one in Brichevo. And so at the age of 10 Boaz left Brichevo and, accompanied by his parents, went to the town Beltz. To be admitted to a state school  one had to pass entrance exams, but Boaz failed.
 It is not clear whether this related to his Jewish origin, or the situation itself caused him to feel too much pressure. Family friends suggested to settle this through their connections, but when this became known to Boaz he began to cry and refused to do things this way. After all, he was admitted to a private school and continued his studies there. During this period, he lived in the family of far relatives.

\begin{figure}[!b]
\vspace*{-2mm}
\centering
\includegraphics[scale=0.57]{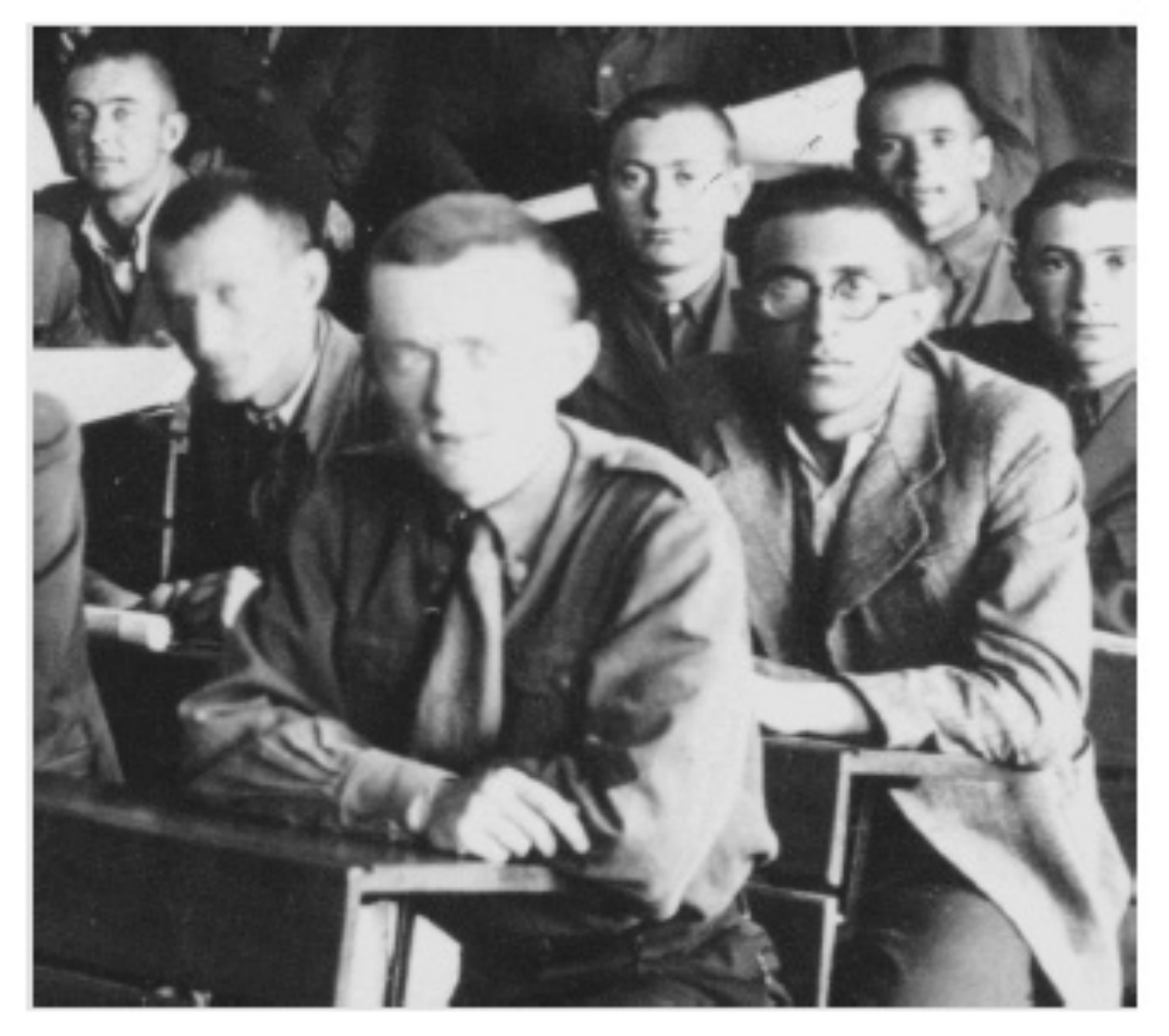}\\
{Boaz in gymnasia ``Tarbut'', 1939}\vspace*{-1mm}
\end{figure}

\medskip
Boaz completed his school education in 1939, in the district's central town Soroca, known for its well-preserved medieval fort, very close to the border with Ukraine. This was a big town with a mainly Jewish population; but Romanians, Germans and Russians also lived there. In Soroca, Boaz studied in the gymnasia ``Tarbut'' (``culture'' in Hebrew) -- part of a network of secular Jewish educational institutions in Eastern Europe established in the period between the two world wars.

Boaz was very successful in his school studies. Besides mathematics, he very much liked history -- general and Jewish. Boaz always emphasized that he was very lucky to have good teachers in both mathematics and humanities. Some of them arrived to Bessarabia from Ukraine, Belarus and Lithuania, with knowledge of different languages.

\begin{figure}[!h]
\centering
\includegraphics[scale=0.34]{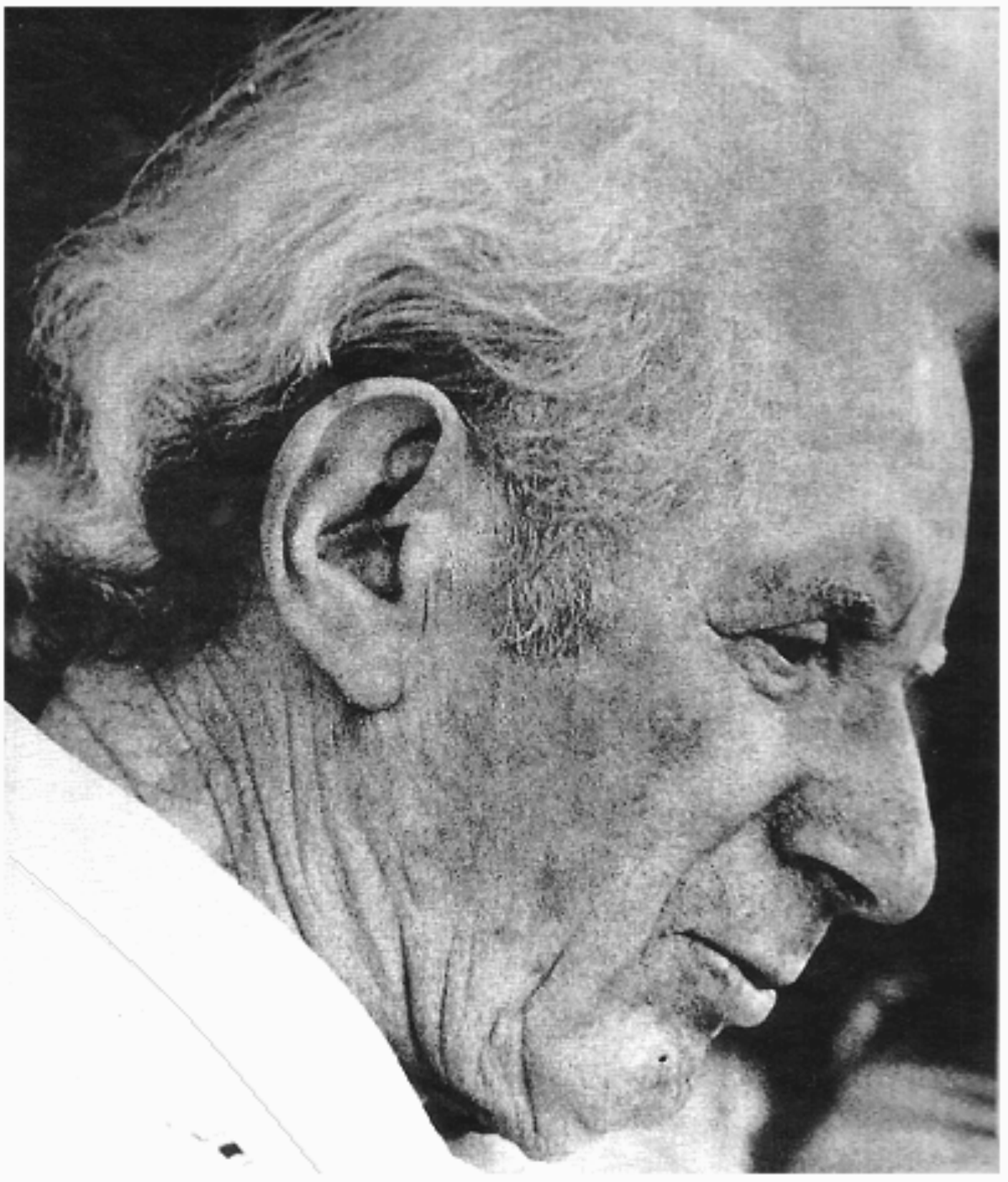}\vspace*{2mm}\\
Aharon Bertini -- Boaz's teacher and director of gymnasia ``Tarbut''
\end{figure}

Besides his first teacher Pinchas Golergant, Boaz especially admired Aharon Bertini. He was one of the brightest examples of those highly educated people that promoted culture and education in Bessarabian towns. Born in Brichevo, he studied in Charles University in Prague and
in Sorbonne, was a director of gymnasia ``Tarbut'', and eventually became a known Israeli writer, poet, translator and educator. He was on editorial boards of several books and encyclopedias, such as History of Brichevo and Encyclopedia of Diaspora (one of its volumes was devoted to Bessarabian Jewry). Boaz shared with the family a variety of recollections about him.

\medskip
In high school, Boaz started to write poems in the Hebrew and Romanian languages; they were even published in some local newspapers. He signed them with the pseudonym Boaz Halakhmi (that means Boaz from Bethlehem), after the name of the great-grandfather of the biblical King David. Interestingly, about 45 years later, the name Halakhmi returned to our family. After repatriation to Israel, Boaz's younger son Yosef decided to take a Hebrew family name, and Boaz suggested that this would be Halakhmi. Few years later, Boaz's nephew (son of his brother) made the same change.

\section{Participation in the Zionist-Socialist movement}

Boaz mentioned that as citizens of Romania, they had the right to vote in parliament elections. His father used to vote for the Agricultural Party, but Boaz himself did not take any part in elections. Rather, he became very active in the Zionist movement that was very popular among young people.

\medskip
There were a variety of groups, and Boaz joined the one that was especially close to his views, called HaShomer HaTzair (Young Guard) -- a Zionist-Socialist movement oriented to humanistic Judaism, Zionism and socialism. The main goals of this movement were repatriation to the Land of Israel, revival of the country and leading a collective life style. Eventually, Boaz even became a member of its Central Committee in Romania. He was involved in organizational and education activities, participated in meetings in various parts of the country, and delivered speeches urging the young people to join and to build socialism in the Land of Israel.

\begin{figure}[!h]
\centering
\includegraphics[scale=0.7]{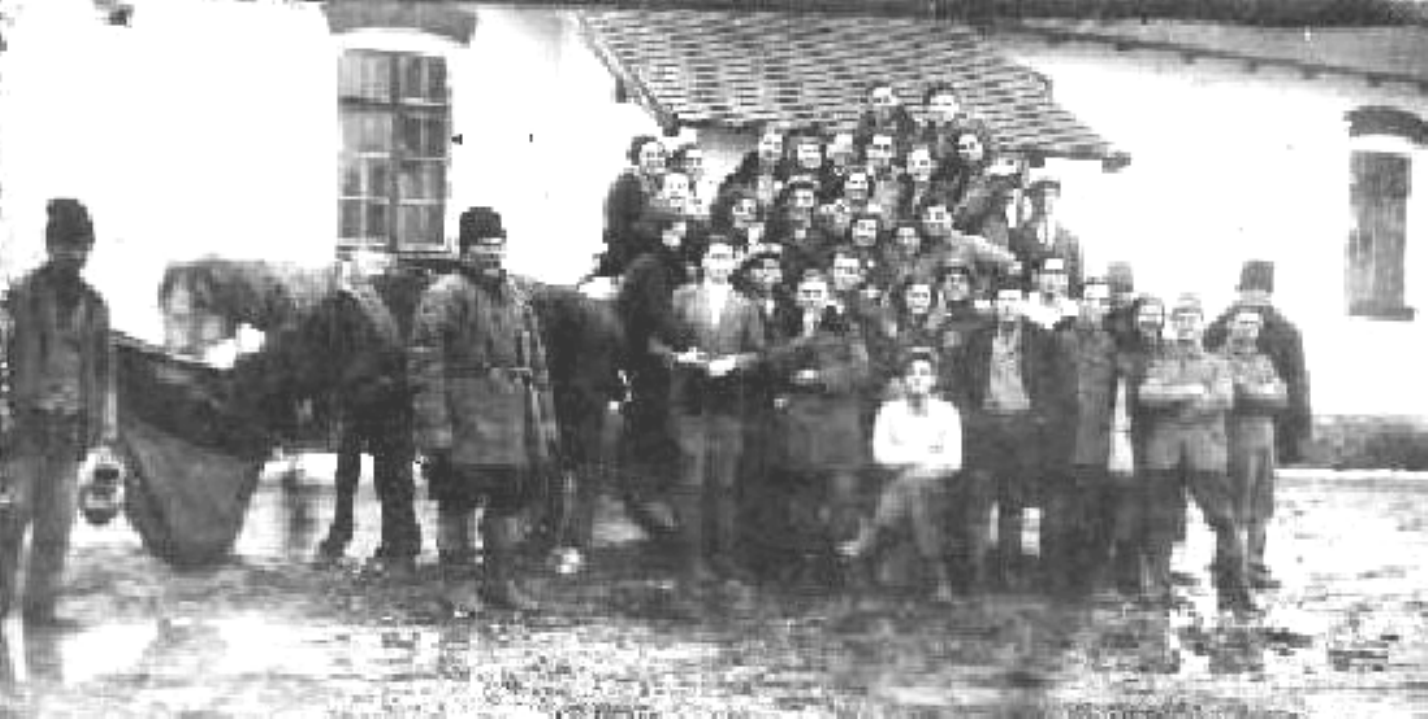}\vspace*{2mm}\\
{Youth of Brichevo -- training before repatriation}
\end{figure}

Boaz also took part in training activities with the goal to develop skills in various types of agricultural work that would prepare him to work in a kibbutz after repatriation. While his father wanted him to study in  the Hebrew University of Jerusalem,  Boaz saw himself as an educated shephard grazing sheep and playing the flute for them, and struggling for equality of all people.

At that time, not many repatriation permissions were given to Jews by British authorities in Palestine. Boaz's estimation was that his turn to get a permission will come somewhere in 1942--43. But very soon dramatic geopolitical changes occurred, and all these dreams became unattainable and even dangerous.

\medskip
Boaz repatriated to Israel only 40 years later. He then put much effort to find his fellow friends from HaShomer Hatzair living in Israel, and was happy to meet them after many years of separation.

\section{Dramatic events before the war}

On August 23 1939, the USSR and Nazi Germany signed a non-aggression pact (known as the Molotov-Ribbentrop pact). According to its secret protocol, the two countries divided between them lands in Eastern Europe. In particular, following this pact and following the ultimatum submitted by the USSR to the Romanian government, on June 28, 1940 the Soviet army entered Bessarabia, and it became a part of the USSR. It was the beginning of a new life.

\medskip
Socialist ideas were popular in Romanian society in the 1930s. Many young people were members of the Communist Party (even though it was outlawed) or of other socialist-oriented parties and movements such as HaShomer HaTzair where Boaz was so active. They wanted socialism here and now, and hence they were quite enthusiastic about the change.

Even when information about political trials and executions in the USSR reached Bessarabia, the idealistic young people like Boaz just refused to believe it. Boaz recalled that when a Jew was accepted by the new authorities to a government job as a postman -- people viewed this as an example of achievements and freedoms that socialism could bring to Jewish population. And he was especially impressed when eradication of illiteracy was announced.

During the first few months after the establishment of the new power in the region, unlimited crossings of borders between Bessarabia and Romania were still possible. However, in Boaz's family the possibility to leave was not considered at all. Boaz's choice was socialism - the hope was that Soviets and HaShomer HaTzair may help each other in achieving their goals. But very soon it became clear that borders became closed, and repatriation would not be possible.

People began to understand the new reality. The new authorities closed Jewish schools in Brichevo; private property was confiscated. In particular, Soviet officials arrived to the shop owned by Boaz's father, and just took for themselves everything they found there; this was done to many locals that owned businesses. People (especially the younger generation) started to leave Brichevo, looking for better opportunities elsewhere.

It was in these circumstances that Boaz moved to Kishinev, where in the fall of 1940 he started his higher education in the Pedagogical Institute, with great hopes that there, nothing would stop him from what he decided to do -- to study mathematics.

\section{Languages}

Besides his major talent in mathematics, Boaz was also very good in mastering various languages. His mother tongue was Yiddish. This was the language used in his family; later Boaz and his wife Berta also used Yiddish when they did not want their children to understand them. During his school years, Boaz also studied Romanian (the state language in Bessarabia till 1940), Hebrew, French, German and even Latin. When he started his higher education in 1940 (after Bessarabia became part of the USSR), he also mastered Russian. During most part of his life, Hebrew and Russian were his main languages.

\medskip
In the USSR, where Hebrew was in fact a forbidden language, Boaz found a way to retain his Hebrew on a high level. First, he used to listen to the Voice of Israel. In general, Soviet authorities used special noise-generation devices to make listening to broadcasts from Western countries extremely difficult. But programs in Hebrew could be listened with no disturbance; probably, it was assumed that nobody knew this language in the USSR. Boaz liked to listen to all kind of programs -- from news to football reportages, even though he never was a football fan. Also, newspapers published by communist parties from all over the world were available in the USSR -- and Boaz used to buy Israel's \textit{Zo Ha-Derech} (``This is the way'') to practice reading in Hebrew and to get some kind of news from Israel.

\medskip
In the 1970s Boaz got from a friend a rare edition of the Bible (published in Russia in pre-Soviet times) with parallel texts in Hebrew and Russian. He suggested to me to open the book on an arbitrary page, hide from him the Russian text (but in a way that I could see it) and to check how he smoothly translates from Hebrew to Russian. His translation was indeed almost identical to canonical Russian text of the Bible. It was especially impressive because biblical language (either Hebrew or Russian) significantly differs from that used in day-to-day life.

When Boaz repatriated to Israel all this definitely helped him. He had no problems in reading newspapers and watching TV, and could freely communicate with people, including his colleagues in Tel-Aviv University. However, for his lectures to students, he had to extend his Hebrew with some specific professional terminology.

As for Russian, Boaz recalled that generally in Brichevo many people were able to say a few phrases in a very basic and clumsy Russian, without even being aware of the mistakes they made. Some people hired private teachers to study Russian; the hope was that this would allow them to go to Russia to find a job there. At that time it meant going abroad!

Things started to change in 1940, and specifically for Boaz when he became a student at the Pedagogical Institute in Kishinev where learning the Russian language was obligatory. The lecturer was a highly cultured person; his Russian was very rich and on a very high level. Boaz made quick progress, and one year later his results in dictations ranked third among all students, including those with Russian roots. But the process of mastering the language was not that simple -- as Boaz recalled, it required a kind of self-control: when he was not sure about the correctness of certain forms of words, he just avoided using them. Eventually Russian became for him another native language.

\medskip
And what about English? Boaz's first exposure to this language happened when he was evacuated during the war; he then was asked to translate to Russian some technical documentation for American equipment. This was a very challenging task for a person who does not know English and is completely ignorant in machinery and relevant terminology. But Boaz's manager promised to help him in coping with technical terminology, and somehow succeeded to convince Boaz that since English is a mixture of French and German (that Boaz had studied), and since he was such a clever guy, he would cope with the task. The strange thing was that they, together, ultimately managed to do the job. As Boaz summarized, this mission proved to be helpful in his future career.

\section{Deportation}

In the night between June 12 and June 13, 1941 a very dramatic event happened: mass deportations of the population from the ``new lands'' that USSR annexed according to Molotov-Ribbentrop pact -- including Bessarabia, and in particular Brichevo. The official reason for this action was ``cleansing'' from the so called ``anti-Soviet, counter-revolutional, criminal and socially dangerous elements''.

\medskip
People were taken out of their homes, transported on horse carts to the nearest train station -- and from there, on trains (prepared by the authorities in advance) they were transferred thousands of kilometers away, to Siberia. Their property was confiscated.

The entire family of Boaz (his parents, sister and brother) was expelled; but he himself avoided the deportation -- he studied in Kishinev and just  was not at home on that date. He became aware of the event on the same day, when he heard rumors about what had happened. He called Brichevo to figure out what was known about his family. As he recalled, in the morning of June 13 all of Bessarabia knew about the deportations and was on phones -- people tried to find out the fate of their relatives and friends. At first, Boaz thought that only wealthy families (like his) were expelled; but very quickly, it became clear that this ``criterion''  was not relevant.

\medskip
According to the laws that existed in the USSR in those times, the families were sent to Siberia for a period from 5 to 20 years. For the family of Boaz, a 10-years term was assigned. Heads of families were separated from their wives and children; they were put on different trains and sent to different places.

This was a very tragic and traumatic event. But in paradoxical and somewhat absurd way, it saved the family from the truly terrible tragedy -- the Holocaust.

\section{During the war}

The Second World War came to the USSR on June 22, 1941 ---  just 9 days after the deportation. Bessarabia was occupied by German and Romanian troops in July 1941. Following the Germans bombardment of Kishinev, the staff and the students of the institute in which Boaz studied, were evacuated from the city far to the East -- first to Orenburg and then to Buguruslan in the Ural Mountains area. There Boaz continued his studies and at the same time worked in a gas company to earn money for a living.

\medskip
And so the family was dispersed over the huge and unknown territory: Boaz, his father Abraham, his mother Brucha with children Ratza and Mendel -- they all were in different places. During the next two years they knew nothing about each other; they did not even know who (and whether) was alive. But they continued their attempts to find each other.

Boaz wrote a letter to Lavrenty Beria --- the minister of internal security and chief of the secret police. Boaz claimed that his father was not a capitalist and was not involved in any political activity; he also mentioned that he himself was a socialist -- and asked for help in finding his family. Eventually Boaz got an answer with his father's address -- a place where he was involved in forced labor in the timber industry. Boaz succeeded to visit his father there, and this was a very emotional meeting; even many years later Boaz could not talk about it calmly. The meeting was short -- in fact, it was forbidden for Boaz to come to this place at all; he did this at his own risk.

In parallel, his mother wrote a letter to the office that held information about all evacuated factories, institutions and people. She could not know of course that this office resided in the same town where Boaz then lived following the evacuation of his institute -- an astonishing coincidence! This way she got Boaz's address.

\medskip
Finally, in 1943, after two years of uncertainty, they could feel some level of confidence about each other's situation -- they knew that at least they were alive, and could exchange letters. Father of Boaz was then allowed to join the rest of his expelled family, not including Boaz -- he was not expelled, and hence he remained separated from them.

\medskip
In 1944 Boaz with his institute returned to Kishinev; the war still continued, but in Europe, beyond the USSR borders. Boaz successfully finished his studies in the institute, and till the end of the war in May 1945 he worked as a math teacher in Beltz -- the town to which he arrived at age 10 for his secondary school studies.

\section{Years in the university}

Boaz accomplished his higher education in mathematics in the years 1945--1947, when he was a student at Chernowitz University. The University was founded in 1875 when Czernowitz was the capital of~the~Duchy~of~Bukovina, {a~Cisleithanian~crown~land} of~Austria-Hungary (in 1940 it became a part of the Western Ukraine area in the USSR). Nowadays, the university is based at~the~Residence~of~Buk\-ovinian and Dalmatian Metropolitans~building complex,
a~UNESCO$\;$World Heritage Site~since~2011.

\begin{figure}[!h]
\vspace{2mm}
\centering
\includegraphics[scale=3.5]{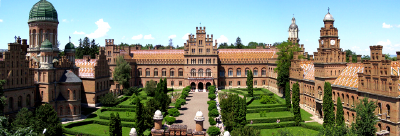}\vspace*{2mm}\\
{Chernowitz University}
\end{figure}

The staff at the faculty of mathematics was mainly with roots in Moscow mathematical research groups, and this guaranteed a high level of teaching. Boaz was especially impressed by Dr. A.Bobrov, a former PhD student of the famous Russian mathematician A.Kolmogorov. He had an open and democratic communication style, with good sense of humor and self-criticism. This was new to Boaz; later he himself adopted this style in his work with students and colleagues.

\medskip
Boaz got a job in the restoration of the very rich university library (during the war, many books and journals were kept in underground rooms). Boaz used this not only as a way to earn some money for a living, but primarily as a great opportunity for serious self-education.

\medskip
But most important was that Bobrov's seminar based on F.Hausdorff's book ``Set Theory'' attracted Boaz, and he decided to do research on descriptive set theory. Bobrov then recommended to him to establish contacts with prominent researchers in Moscow -- A.Kolmogorov, P.Novikov and A.Lyapunov; first meetings with them took place in 1946--47. Here it is important to emphasize again the great role Novikov and Lyapunov played in Boaz's professional life, and their high human and moral standards. Boaz recalled again and again the warm welcome he received at the homes of these scientists and their uncompromising support in difficult periods of his life.

\begin{figure}[!h]
\centering
\includegraphics[scale=0.28]{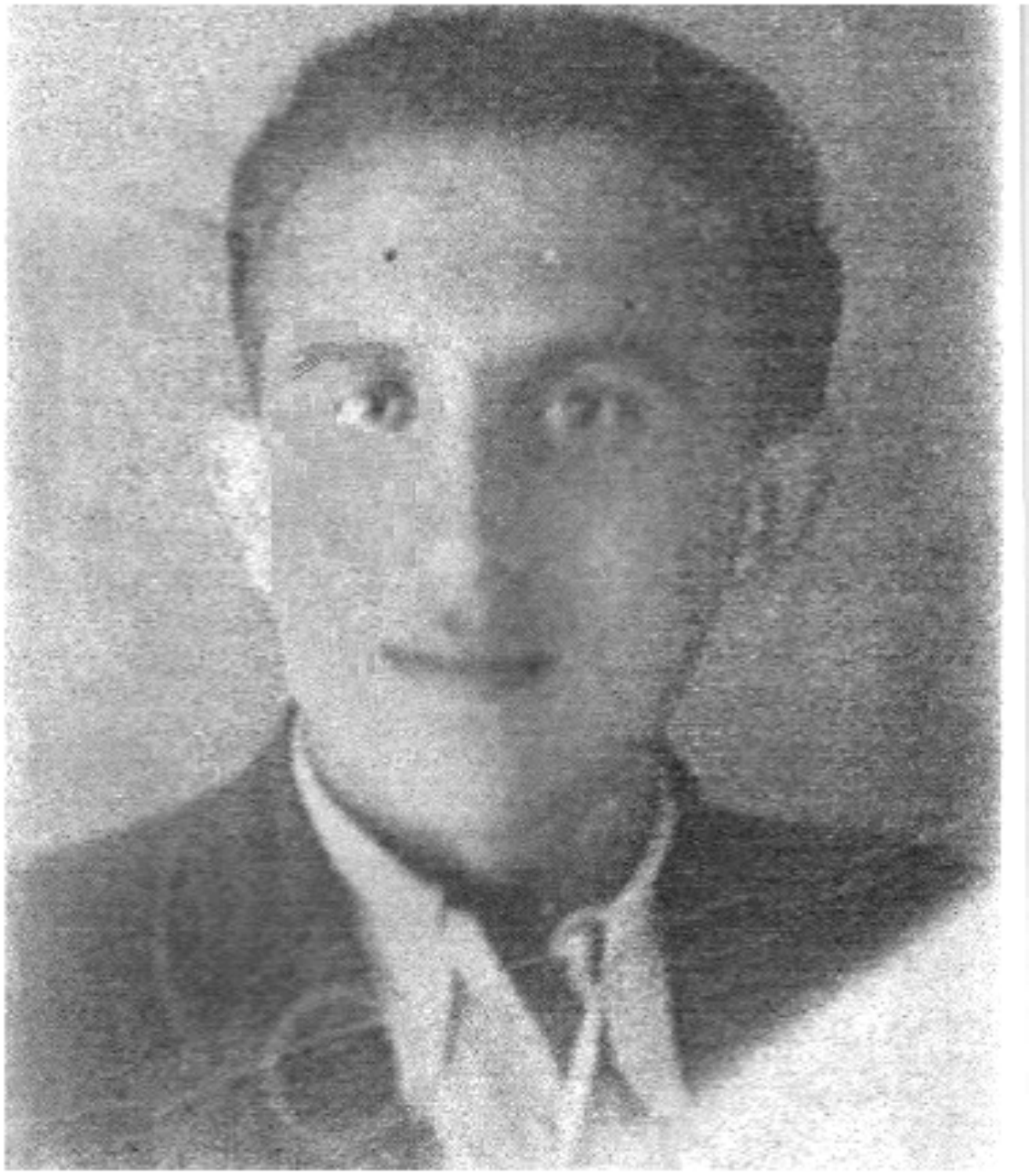}\\
Boaz -- student at Chernowitz University
\end{figure}

\begin{figure}[!h]
\vspace*{-2mm}
\centering
\begin{tabular}{c}
\includegraphics[scale=0.6]{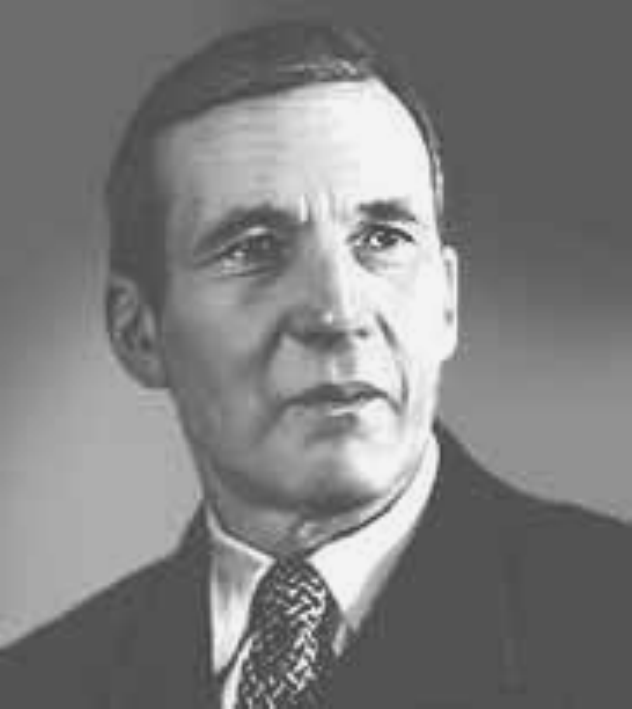} \hspace{1cm} \includegraphics[scale=0.49]{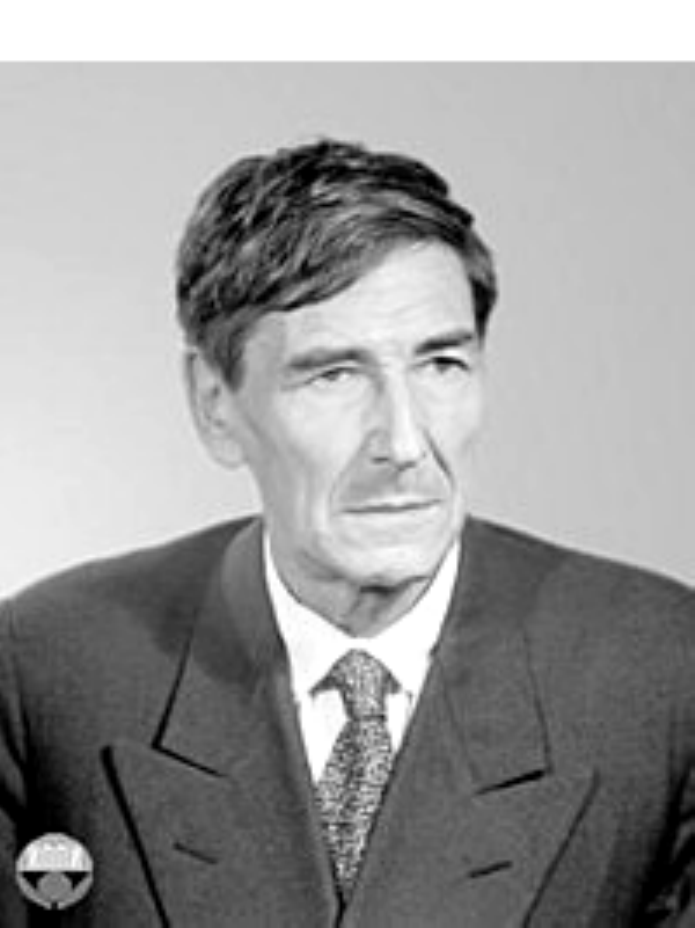} \hspace{1cm} {\raisebox{-1ex}{\includegraphics[scale=0.88]{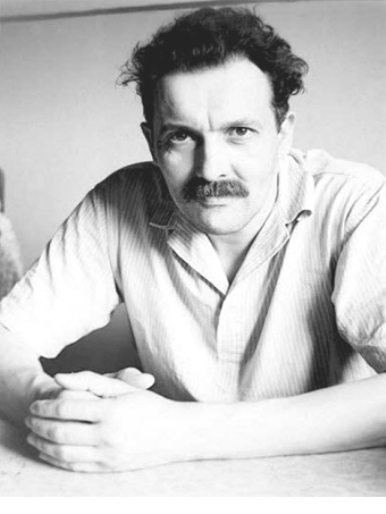}}}\\
\hspace*{-0.3cm} Andrei N. Kolmogorov  \hspace*{1.8cm}\hfil    Petr S. Novikov \hspace*{1.1cm} \hfil   Alexey A. Lyapunov
\end{tabular}
\label{fig:discussion}\vspace*{-3mm}
\end{figure}

\medskip
This is how Boaz described his first meetings with Lyapunov and Novikov: ``Imagine who I was at that time: a Jewish student, a stranger who in fact is quite new in this country and whose family was deported to Siberia (even though they didn't know about this)$\ldots$  And they accepted me in a very friendly way, with respect to my first scientific interests, and gave me good advices and suggested new research areas. This was far from being obvious at that time in the USSR''.

\begin{figure}[!b]
\vspace{-2mm}
\centering
\includegraphics[scale=0.35]{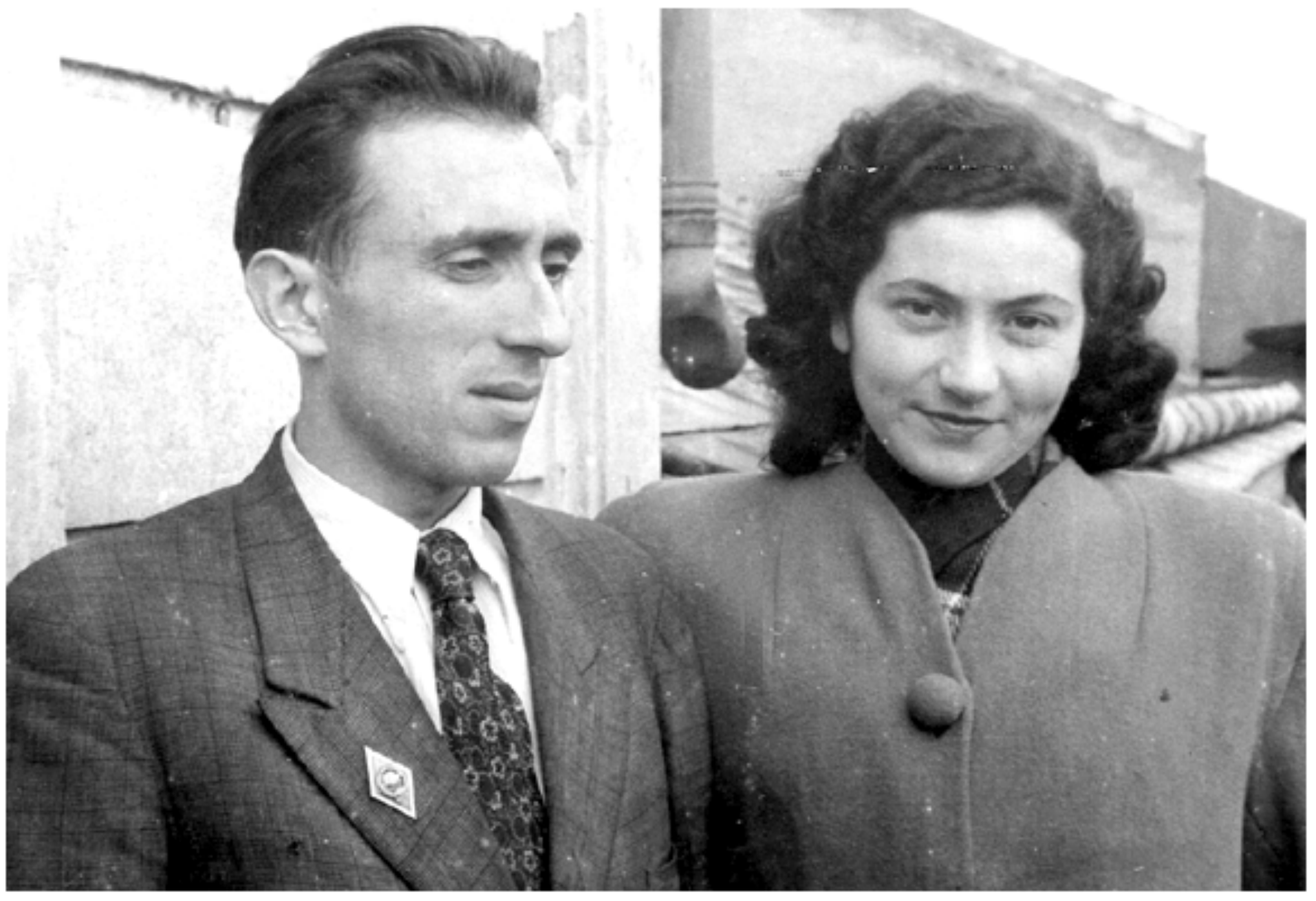}\vspace*{1mm}\\
{Boaz and Berta in 1947}\vspace*{-1mm}
\end{figure}

After his short meeting with A.Kolmogorov in Moscow in 1946, Boaz continued to Siberia. Officially, he was sent by Chernowitz University to look for potential math students in the East of the country. But for him the main goal was to bring the family from exile back to Chernowitz -- Boaz's father was freed after 5 years, with mother, sister and brother allowed to join him.

For Boaz, there were three reasons to feel that 1946 was a happy year. First, the family survived and was freed. Second, he started his first research activities, and established connections with leading scientists in Moscow. And third -- he met Berta Rabinovich, his future wife.

\section{PhD student in Kiev and Moscow}

In May 1947, shortly after he got married with Berta, Boaz met with P.S.Novikov at his home in Moscow. It was at this meeting that Novikov attracted Boaz's attention to new research areas: logic and computability theory. Novikov also suggested to Boaz to be his scientific advisor. Boaz gladly accepted this, and in October 1947 he became a PhD student in Kiev, in the Institute of Mathematics of the Ukrainian Academy of Science. Interestingly, at that time the director of the Institute was Michael A. Lavrentyev, who ten years later became one of the founders and the first President of the Siberian Branch of the USSR Academy of Sciences in Novosibirsk -- where Boaz worked from 1960 till his repatriation to Israel at the end of 1980.

There was a curious episode related to the qualifying exam that Boaz had to pass as a PhD candidate. His examiners
  (M.Lavrentyev, N.Bogolyubov, S.Krein) were involved in research in applied mathematics, and so for the exam Boaz decided to choose topics related to their areas of interest. But he failed in his first attempt to pass the exam! After at age 10 he had failed in the entrance exam to a secondary school in Beltz, this was his second (and the last) failure. He undertook one more attempt, and passed with excellence.

\begin{figure}[!h]
\centering
\includegraphics[scale=1.25]{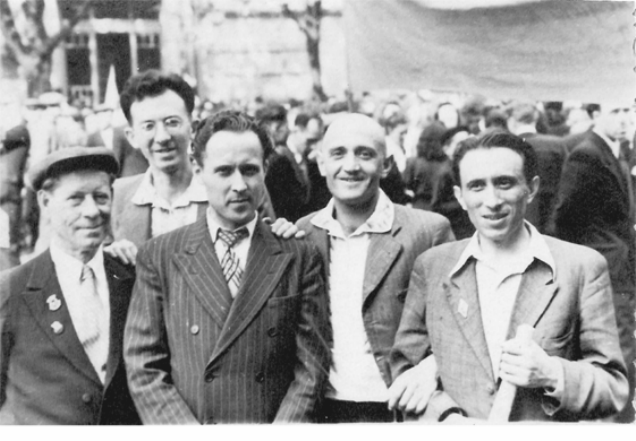}\vspace*{2mm}\\
{Boaz (right) -- a PhD student in Kiev, with colleagues}
\end{figure}

After a letter written by Novikov to Lavrentyev, specialization of Boaz in logic was approved, including long-term visits from Kiev to Novikov in Moscow. In Moscow he participated in the so called Big Seminar headed by P.Novikov and S.Yanovskaya. Boaz presented there the results of his future PhD thesis, and S.Yanovskaya suggested to provide an official support letter when his thesis would be ready.

\medskip
The seminar participants had discussed a wide spectrum of subjects, including discussions related to the books on logic by Tarski and by Hilbert and Accerman translated to Russian in 1947--48. The official academic establishment was hostile to mathematical logic, and in particular, in the preface to the translation of Tarski's book he was blamed to be a militant bourgeois. This is important for understanding the atmosphere related to Boaz's area of research. In particular, couple of years later, when in 1950 he started to work in Penza after defending his PhD thesis, Boaz was accused of being ``a Carnap-type idealist''. Only the strong support that he got from his prominent colleagues in Moscow helped him to overcome that dangerous situation; he used to say ``they saved me''.

\begin{figure}[!h]
\centering
\includegraphics[scale=0.39]{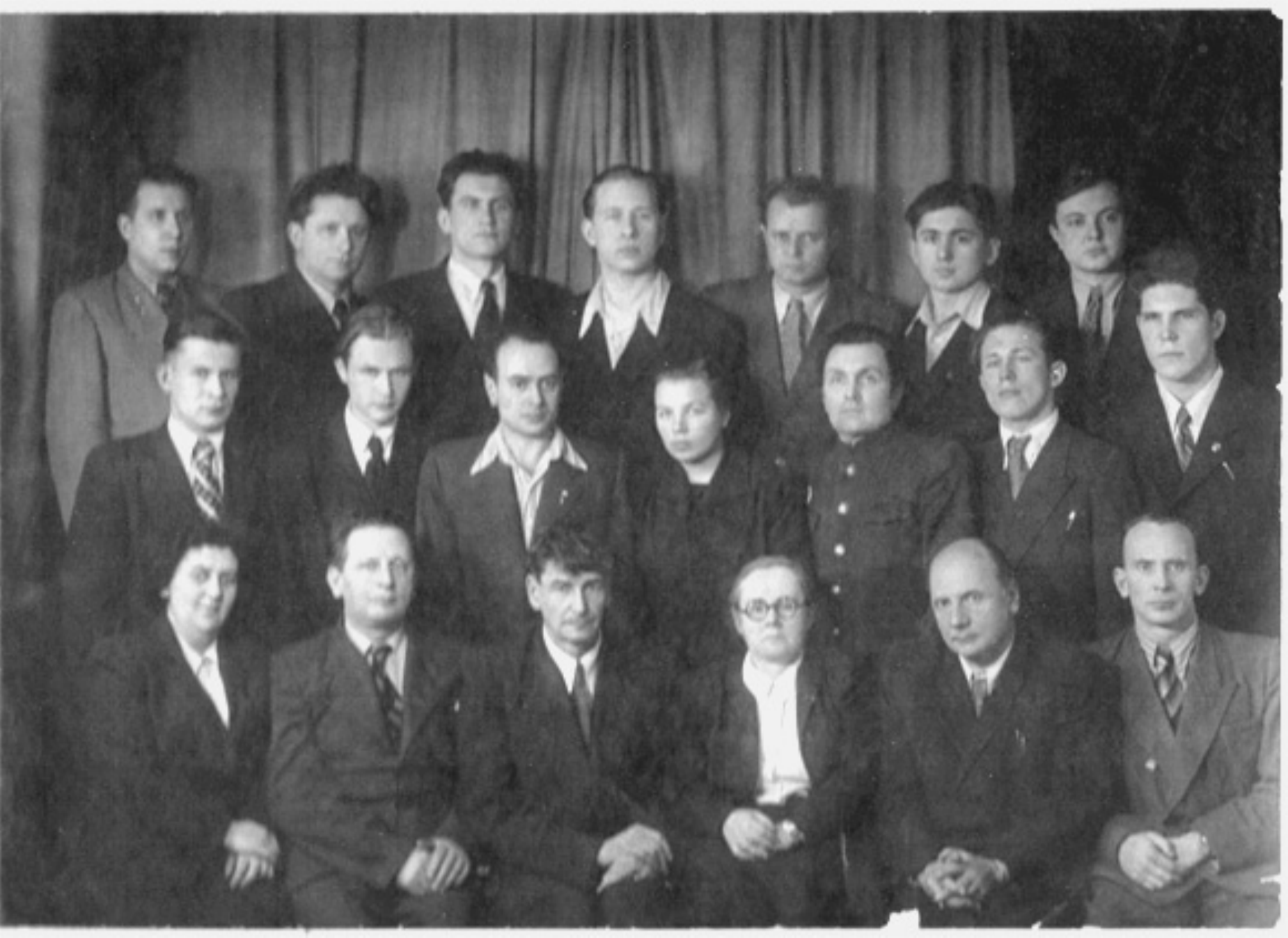}\vspace*{1.5mm}\\
{Big Seminar in Moscow. First row in center: P.Novikov and S.Yanovskaya}
\end{figure}

\medskip
With some of the participants of the Big Seminar, Boaz kept fruitful and friendly relations during many years. Among them -- his close colleagues Alexander Kuznetsov and Alexei Gladky; the latter eventually (in 1967) joined Boaz in the ``Automata theory and mathematical linguistic'' department that Boaz headed in the Institute of Mathematics in Novosibirsk. With another one Boaz got familiar even earlier, in 1947 in Chernowitz University. This was Alexander Esenin-Volpin, son of the famous Russian poet Sergei Esenin. He arrived to Chernowitz after he got his PhD in mathematical logic in Moscow; but even before he managed to work a single day, he was arrested and sent to exile. Over the years, he was not only a mathematician, but also
a notable dissident, political prisoner and human rights fighter. Boaz's contacts with him (sometimes direct and sometimes through his mother, a poet Nadezhda Volpin) continued for many years, and their last meeting was in Jerusalem in 2009.

\begin{figure}[!h]
\centering
\includegraphics[scale=1.3]{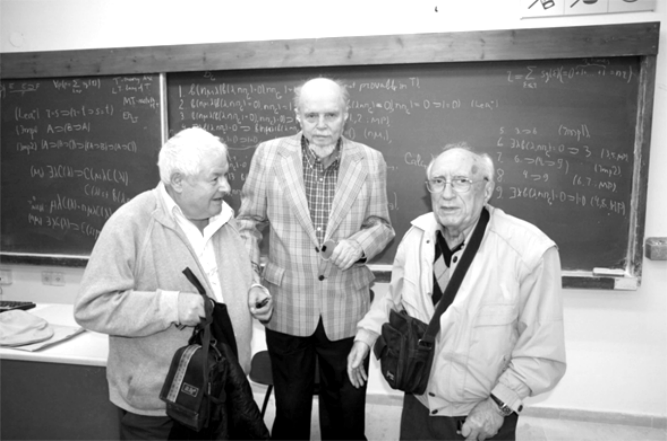}\vspace*{1.5mm}\\
{Boaz (right) and A. Esenin-Volpin (center) in Jerusalem, 2009}\vspace*{-2mm}
\end{figure}

\section{PhD dissertation}

Boaz finished his work on the PhD thesis in 1950. He achieved very good results; the main one is known as Trakhtenbrot's Theorem:  \textit{The validity of first-order statements that hold true for all finite universes is undecidable}. The date for the thesis defense in Kiev was already set, and all the needed preparations have been completed.

\medskip
However, very close to the planned defense, the fate prepared a new test for Boaz -- his family was expelled again. Somehow, the authorities decided that it was a mistake to let the family return in 1946 to Chernowitz, and that they should go to Siberia again. By that time, Boaz's sister had already married, with a baby boy -- but this did not change that striking decision. Authorities were also looking for Boaz, in order to send him to Siberia as well (he himself was already a father). And it was exactly as it happened in 1940 when the deportations took place -- he was just not with the family at that moment, and hence he again avoided this fate. The information about what happened reached Boaz during his meeting with Novikov at his home. Boaz could not share this information with anyone; hence, he told Novikov that his wife Berta was ill and that he immediately had to go to Chernowitz.

This event was a serious blow for Boaz -- and not only with respect to his family. He was afraid that if things became known in his institute then he could forget about his PhD and his research career would be destroyed. The uncertainty was very high, and the stress was hard to bear. But he was lucky again -- no one knew what happened, and the defense took place as planned. The opponents  (A.Kolmogorov, A.Lyapunov, B.Gnedenko)
 provided excellent reviews that praised his achievements. Boaz summarized this story succinctly: ``I succeeded to slip through''.

\section{Start of work in Penza}

After getting his PhD degree, Boaz started to look for a position in various institutions. It was an extremely bad time for work searches. Starting from 1948, an aggressive anti-Jewish campaign (called ``the struggle against cosmopolitism'') was launched in the USSR. It started with the murder of the great actor Solomon Michoels, and did not stop there -- many Jewish writers, journalists, actors, doctors etc. were arrested and jailed, or even killed. For a Jewish person with an academic background, getting a job in this period (that continued till Stalin's death in March 1953) was next to impossible, definitely in more central areas of the country.

\begin{figure}[!h]
\vspace*{1mm}
\centering
\includegraphics[scale=0.42]{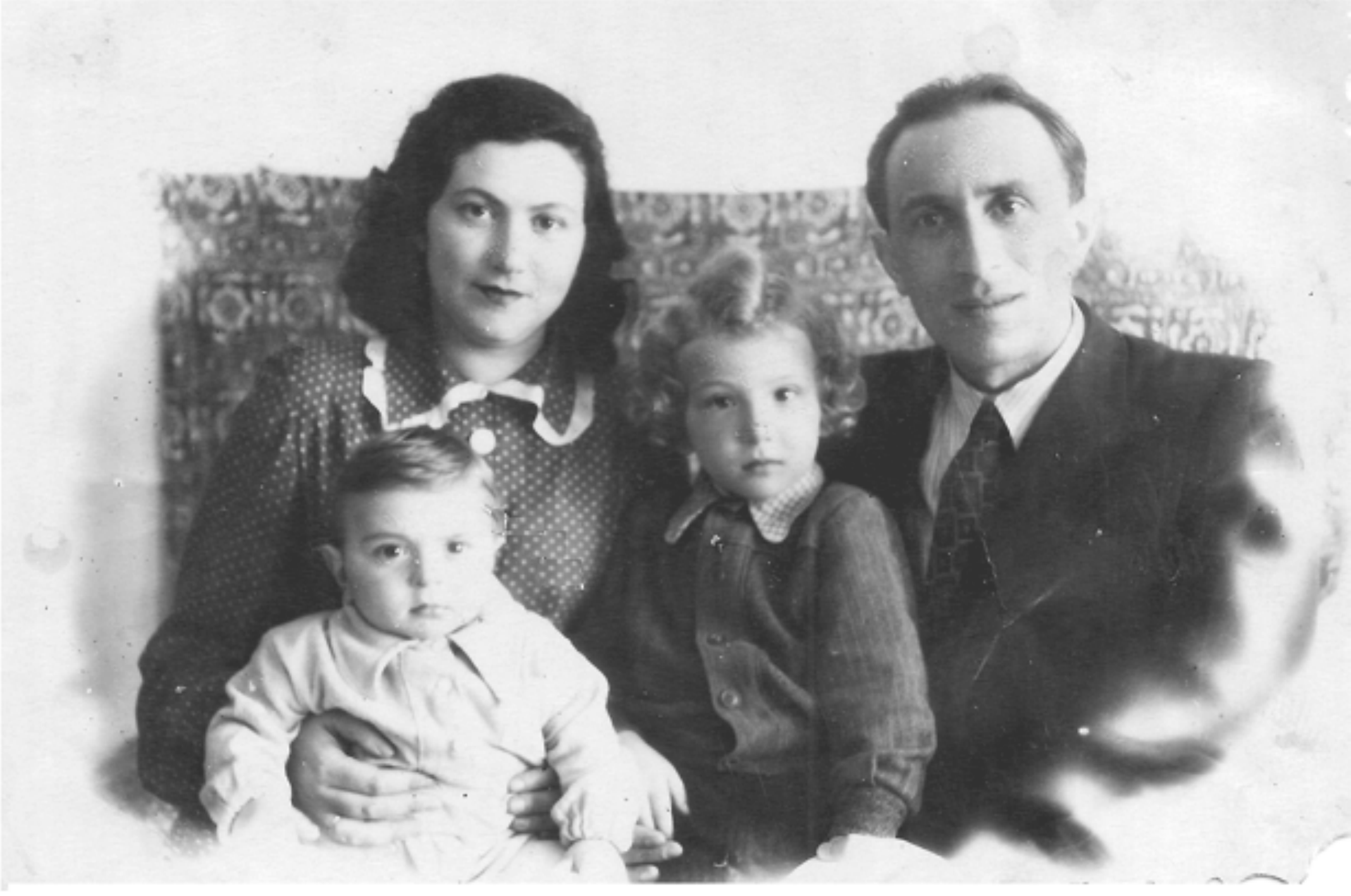}\vspace*{2mm}\\
{Boaz with his family in Penza in early 1950s}\vspace*{-4mm}
\end{figure}

\medskip
Boaz quickly understood that there is no chance in Kiev and in Ukraine. Following the advice by A.Lyapunov, he concentrated his search efforts in the European part of Russia. He sent more than 60 applications to various high education institutions; he got few negative responses, but
 in most cases, there were no answers at all.
Suddenly he got a positive response from the Pedagogical Institute in Kostroma (340 kilometers from Moscow); he arrived in
this town for a meeting with the management of the Institute, but on the next day he was told that actually there was no vacancy.

The situation seemed to be hopeless. Lyapunov then contacted the deputy minister of higher education whom he knew personally, and asked to help. This way Boaz got an offer for a position in the Pedagogical Institute in Penza (650 kilometers from Moscow), with a promise to provide an apartment for his young family. He arrived in Penza on December 5, 1950; his wife Berta with their baby son joined in April 1950. However, the Institute could not provide them a separate apartment; instead, they got two rooms in an apartment shared with another family. This was a very modest accommodation: they had no indoor facilities and needed to burn firewood or coal for heating. But after all, at their age 30, the young family could finally settle down.
	
\vspace{2.5mm}
Ten years in Penza were very fruitful for Boaz. He conducted there his research in logic and automata theory; he published his highly rated book
``Algorithms and automatic computing machines'' (in Russian:\! {\selectlanguage{russian}\textcolor{black}{“Алгоритмы и машинное решение задач”}}) which was very popular, and was translated to a dozen of different languages; and he became known to the research community in the USSR and abroad.

\begin{figure}[!h]
\vspace*{-3mm}
\centering
\includegraphics[scale=0.43]{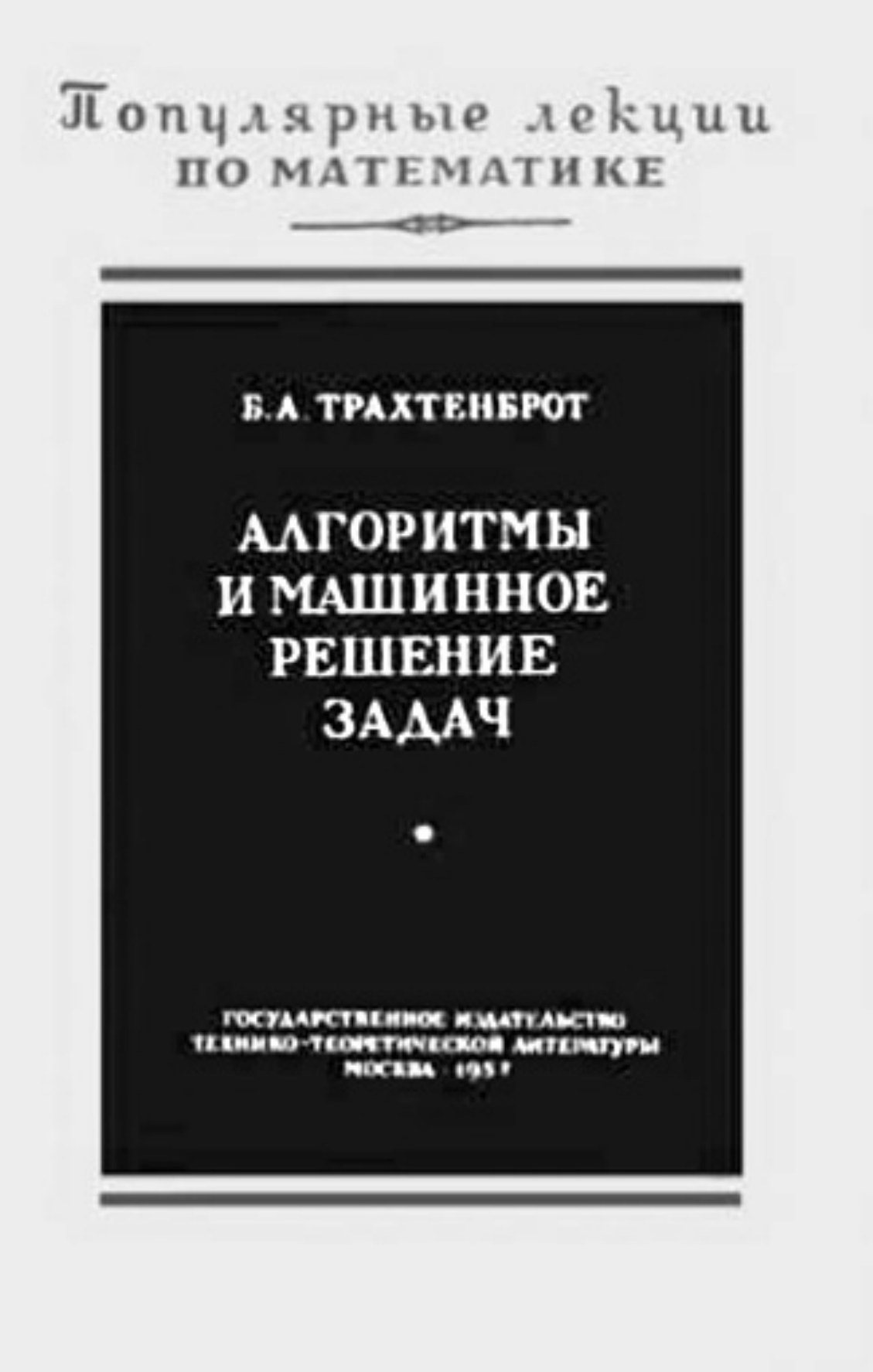} \hspace{9mm} {{\includegraphics[scale=0.309]{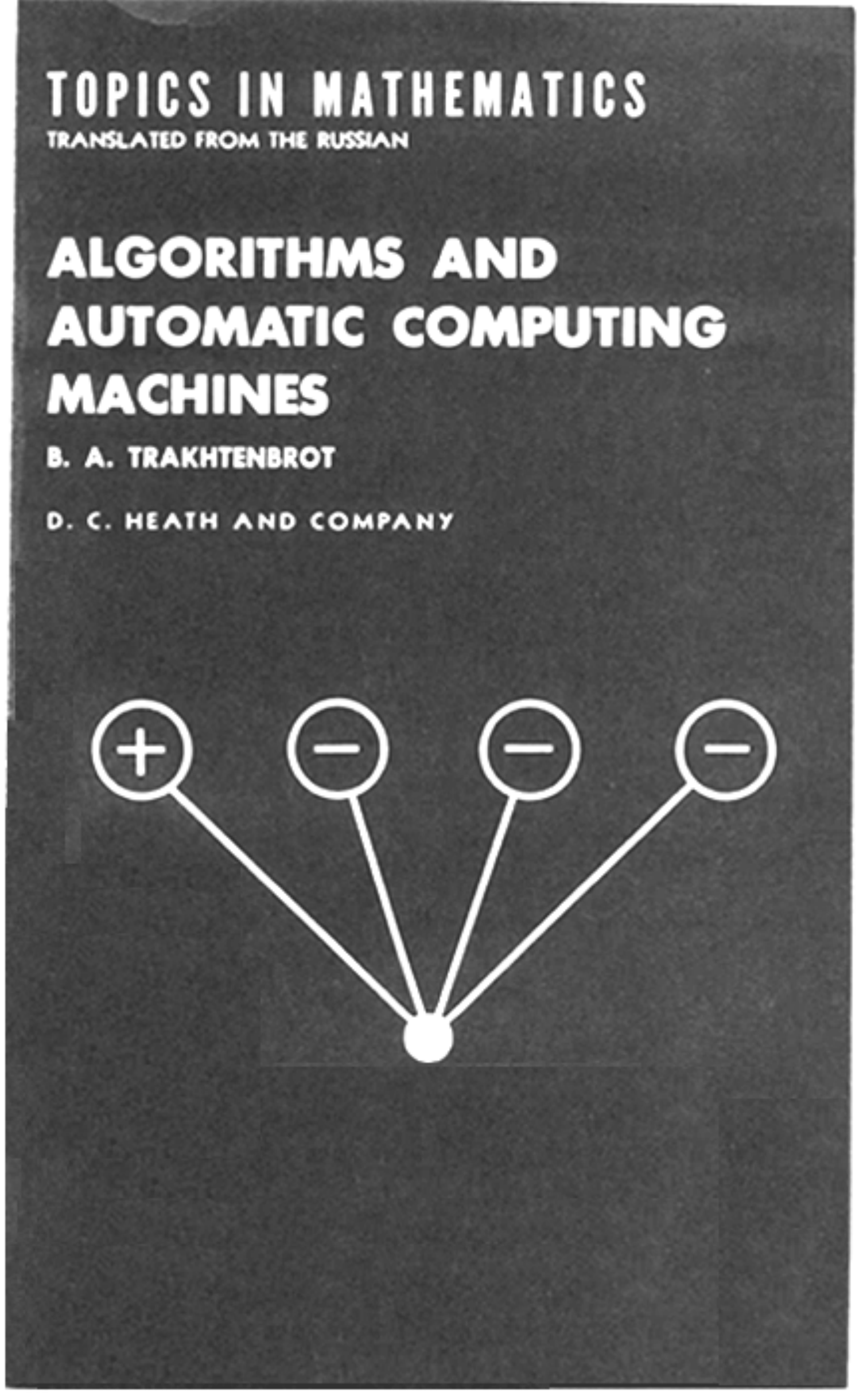}}}\vspace*{-2mm}
\label{fig:discussion}\vspace*{-2mm}
\end{figure}

\eject
And although there were more difficult times and challenges ahead (like already mentioned accusations of idealism, and not only), this was the start of a fascinating road that led him to his great scientific achievements and to general recognition as one of the founding fathers of Theoretical Computer Science.

\section{Chronology of Boaz's life}

\begin{description}
\itemsep=-1pt
\item[] 1921 -- Born in settlement Brichevo, in Bessarabia (Romania)

\item[] 1926-1930 -- Studied first in a \textit{cheder}, and then in an elementary school.

\item[] 1931-1940 -- Studied at a secondary school in Beltz and then in gymnasia ``Tarbut'' (high school) in Soroca

\item[] 1938-40 -- Active member of the youth Zionist-Socialist movement HaShomer HaTzair (Young Guard) in Romania. Preparations for repatriation to Palestine (that was never achieved then).

\item[] 1940, June 28 -- Bessarabia is annexed to the USSR; new life began in a new country, with a new language.

\item[] 1940 -- Began higher education in Kishinev, at the Pedagogical Institute.

\item[] 1941, June 13 -- The family (parents, sister and brother) was deported from Brichevo to Siberia.

\item[] 1941-1944 -- After the German attack on the USSR, evacuated together with the Institute to Buguruslan in the Ural Mountains area. In 1943 succeeded to find his deported family.

\item[] 1944-1945 -- Returned from evacuation and worked as a math teacher in Beltz

\item[] 1945-1947 -- Completed his higher education at Chernowitz University. The family allowed to return from Siberia.

\item[] 1947-1950 -- PhD student in Kiev (Ukraine), with regular visits to his scientific advisor P.S.Novikov in Moscow.

\item[] 1950 --  PhD defense. The main result is known as Trakhtenbrot's Theorem.\\
Parents, sister and brother again deported to Siberia, this time from Chernowitz.\\
Difficulties in finding a job in Ukraine and Russia, due to the government-led anti-Semitism policy (under the slogan of ``struggle against cosmopolitism'').

\item[] 1950-60 -- Worked in Penza as an Associate Professor, first at the Pedagogical Institute and then at the Polytechnic Institute.

\item[] 1951 -- Accused by the management of the Higher Mathematics department as an idealist of the Carnap type, following his lecture ``The Method of Symbolic Calculus in Mathematics'' -- an extremely dangerous (not only for work) accusation in the era of Stalin's paranoia.

\item[] 1960-1980 -- Professor at the Institute of Mathematics of the Siberian Branch of the USSR Academy of Science, in the Novosibirsk Academgorodok. Doctoral dissertation (1962). Professor at Novosibirsk State University.

\item[] 1967-1977 -- Head of the Department of Automata Theory and Mathematical. Linguistics at the Institute of Mathematics.

\item[] 1980, December 26 -- Repatriation to Israel.

\item[] 1981-1991 -- Professor, School of Computer Science at Tel Aviv University.

\item[] 1991-2016 -- Professor emeritus.
\end{description}

\section{Anniversaries and special events}

\begin{description}
\itemsep=-1pt
\item[] 1971 -- Conference at the Institute of Mathematics of the Siberian Branch of the USSR Academy of Science, on the occasion of Boaz's 50th anniversary.

\item[] 1981 -- All-Israel seminar in honor of Boaz's 60th anniversary and his repatriation to Israel.

\item[] 1991 -- International conference at Tel Aviv University in honor of Boaz's 70th anniversary, with participation of leading scientists in the field of Computer Science.

\item[] 1997 -- International colloquium and honorary doctorate (Doctor Honoris Causa) from the University of Jena, Germany -- celebrating Boaz's 75th anniversary.

\item[] 2001 -- In honor of the 80th anniversary, Boaz was invited as a Keynote Speaker at a joint meeting of two leading conferences in the field of computer science theory: ICALP'01 and STOC'01 (Crete, Greece).

\begin{figure}[!h]
\vspace*{2mm}
\centering
\includegraphics[scale=0.55]{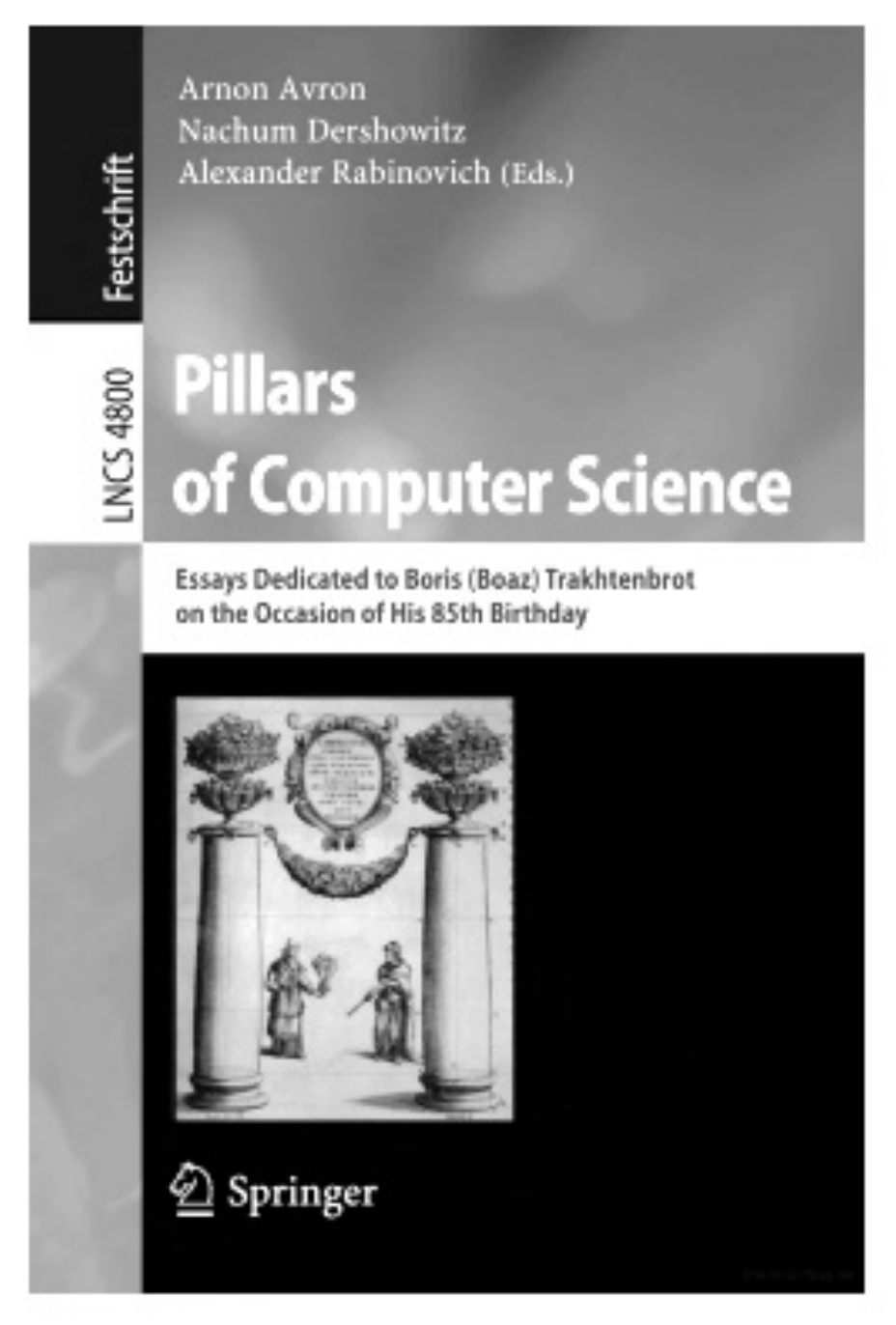} \hspace{9mm} \includegraphics[scale=0.396]{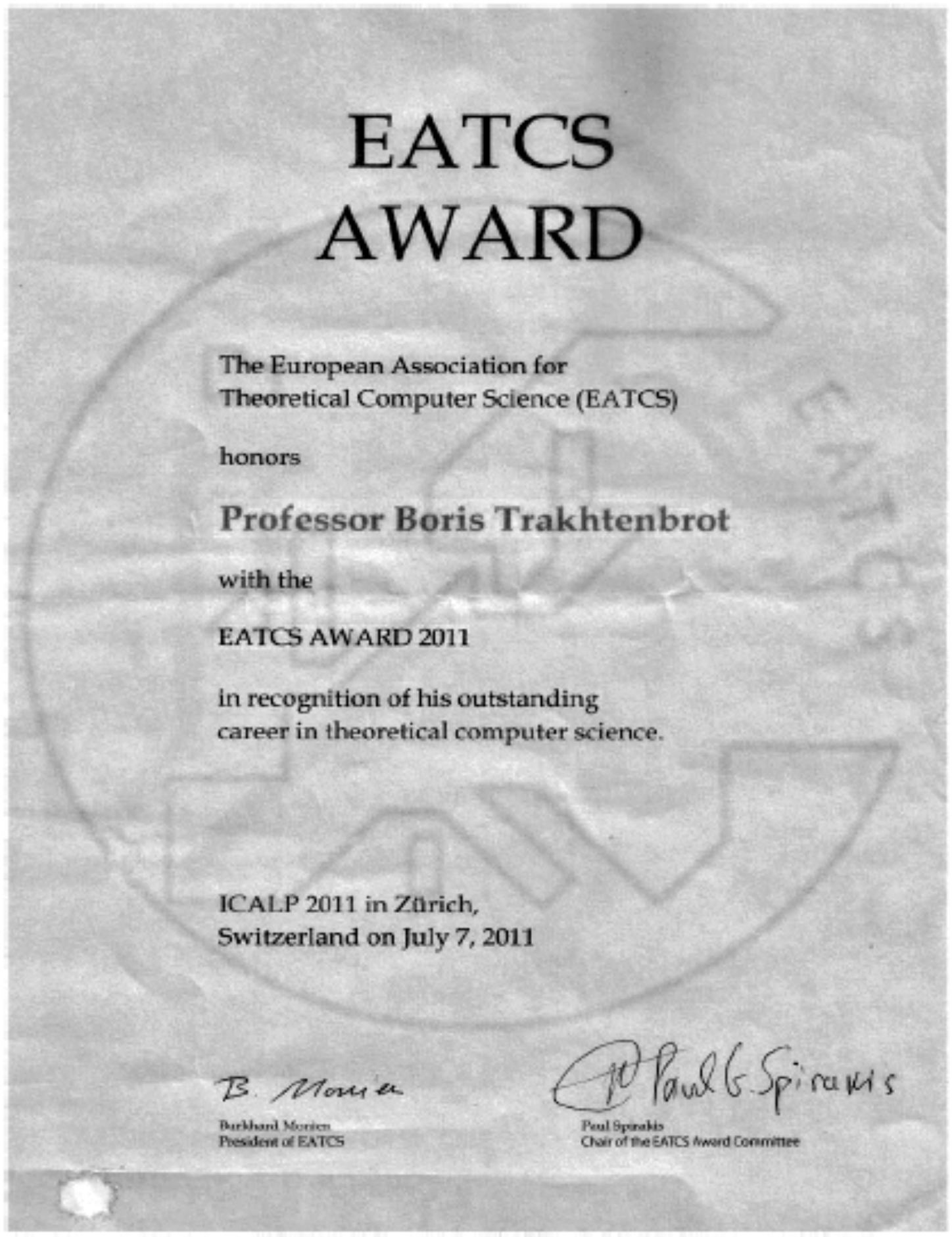}
\label{fig:discussion}
\end{figure}

\item[] 2006 -- Computation Day in honor of Boaz's 85th anniversary at Tel Aviv University. A volume published on this occasion in the Pillars of Computer Science series, which includes both historical reviews and scientific articles by leading researchers -- colleagues and former students of Boaz.

\item[] 2011 -- EATCS Award in recognition of his outstanding carrier in theoretical computer science; presented at the ICALP'11 conference in Zurich.
\end{description}

\end{document}